\documentclass[12pt]{article}
\usepackage{graphicx}
\usepackage{amsfonts,amsmath,amssymb,amsthm}

\begin{document}

\theoremstyle{plain}
\newtheorem{theor}{Theorem}
\newtheorem{lemma}{Lemma}

\theoremstyle{definition}
\newtheorem{defn}{Definition}

\theoremstyle{remark}
\newtheorem*{rem}{Remark}

\def\Z{{\mathbb Z}}
\def\R{{\mathbb R}}
\def\C{{\mathbb C}}

\title{Tetromino tilings and the Tutte polynomial}

\author{
  {\small Jesper Lykke Jacobsen${}^{1,2}$} \\[1mm]
  {\small\it ${}^1$Laboratoire de Physique Th\'eorique
  et Mod\`eles Statistiques}                             \\[-0.2cm]
  {\small\it Universit\'e Paris-Sud, B\^at.~100,
             91405 Orsay, France}                        \\[1mm]
  {\small\it ${}^2$Service de Physique Th\'eorique}      \\[-0.2cm]
  {\small\it CEA Saclay, Orme des Merisiers,
             91191 Gif-sur-Yvette, France}               \\[-0.2cm]
  {\protect\makebox[5in]{\quad}}  
  \\
}

\maketitle
\thispagestyle{empty}   

\begin{abstract}

  We consider tiling rectangles of size $4m \times 4n$ by T-shaped
  tetrominoes.  Each tile is assigned a weight that depends on its
  orientation and position on the lattice. For a particular choice of
  the weights, the generating function of tilings is shown to be the
  evaluation of the multivariate Tutte polynomial $Z_G(Q,{\bf v})$
  (known also to physicists as the partition function of the $Q$-state
  Potts model) on an $(m-1) \times (n-1)$ rectangle $G$, where the parameter
  $Q$ and the edge weights ${\bf v}$ can take arbitrary values depending on
  the tile weights.

\end{abstract}

The problem of evaluating the number of ways a given planar domain can be
covered with a prescribed set of tiles (without leaving any holes or
overlaps) has a long history in recreational mathematics, enumerative
combinatorics, and theoretical physics.

\begin{figure}
  \centering
  \includegraphics[width=150pt]{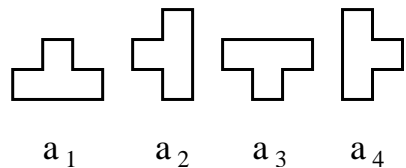}
  \caption{Weights of the four types of T-shaped tetrominoes.}
  \label{fig1}
\end{figure}

\begin{figure}
  \centering
  \includegraphics[width=120pt]{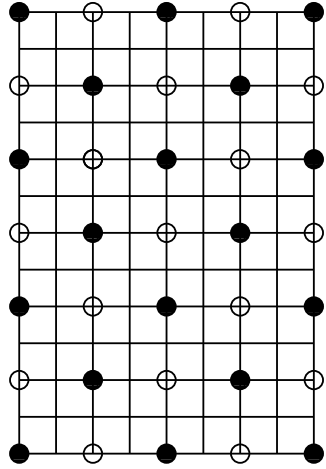}
  \caption{The domain ${\cal D}$ to be tiled (here with $m=2$, $n=3$). Some of
the vertices are distinguished by black or white circles.}
  \label{fig2}
\end{figure}

In this note we consider tiling an $M \times N$ rectangle of the
square lattice (henceforth denoted ${\cal D}$) with T-shaped
tetrominoes, i.e., tiles of size four lattice faces in the shape of
the letter T. The four different orientations of the tiles are shown
in Fig.~\ref{fig1}. Note that the boundary of each tile comprises ten
vertices, of which eight are corners.  Fig.~\ref{fig2} shows the
domain ${\cal D}$ to be tiled. The circles distinguish two subsets of the
vertices which we shall henceforth refer to as either black (filled
circles) or white (empty circles).

We have the following remarkable

\begin{theor}[Walkup \cite{Walkup}]
${\cal D}$ is tileable iff $M=4m$ and $N=4n$. In a valid tiling,
no tile corner covers a white vertex, whereas black vertices are covered
by tile corners only.
\end{theor}

\begin{figure}
  \centering
  \includegraphics[width=75pt]{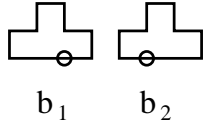}
  \caption{Weights depending on the relative position of a white vertex
within a tile.}
  \label{fig3}
\end{figure}

Theorem 1 implies that tiles can be distinguished not only by their
orientation (cf.~Fig.~\ref{fig1}) but also by their position relative
the the white vertices. There are two possible situations, as shown in
Fig.~\ref{fig3}, depending on which of the two cornerless vertices of
the tile covers the white vertex. (This applies to any of the four
orientations of the tiles, although Fig.~\ref{fig3} only illustrates
the first orientation.)

Note that the lattice $B$ of black vertices is a tilted square lattice of
edge length $2\sqrt{2}$. This can be divided into two sublattices,
$B_{\rm even}$ and $B_{\rm odd}$, which are both straight square lattices
of edge length $4$. Our convention is that all black vertices on the
boundary of ${\cal D}$ belong to $B_{\rm odd}$.

We similarly divide the lattice $W$ of white vertices into two sublattices
$W_{\rm even}$ and $W_{\rm odd}$. Referring again to Fig.~\ref{fig2}, we
adopt the convention that white vertices on the lower and upper boundaries of
${\cal D}$ belong to $W_{\rm odd}$ and that those on the left and right
boundaries belong to $W_{\rm even}$.

Equivalently, any vertex in $B$ or $W$ belongs to the odd (resp.\ even)
sublattice whenever  its distance from the bottom of ${\cal D}$ is twice
an even (resp.\ odd) integer.

Finally, we define the generating function
\begin{equation}
 F_{\cal D}(\{a\},\{b\}) = \sum_{T \in {\cal T}({\cal D})}
 \prod_{t \in T} a_{o(t)}^{w(t)} \prod_{j=1}^2 b_j^{B_j}
 \label{genF}
\end{equation}
where ${\cal T}({\cal D})$ is the set of all T-tetromino tilings
of ${\cal D}$, and $t$ is a single tile within the tiling $T$.
The weights $a_{o(t)}^{w(t)}$ depend on the orientation $o(t)=1,2,3,4$
of the tile $t$ (cf.~Fig.~\ref{fig1}) and on the (unique) white vertex
$w(t)$ that it covers. Moreover, $B_j$ is the number of
tiles of position $j=1,2$ relative to the white vertices (cf.~Fig.~\ref{fig3}).
Note that the weights $b_j$ are not vertex dependent.

\begin{figure}
  \centering
  \includegraphics[width=75pt]{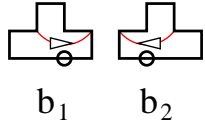}
  \caption{Two ways of decorating a tile by an oriented half-arch, depending on
its position relative to a white vertex.}
  \label{fig4}
\end{figure}

Evidently, the tiles covering the white vertices of the boundary ${\cal B}$
of the domain ${\cal D}$ have their orientation fixed (their long side is
parallel to ${\cal B}$). It follows that those boundary tiles contribute the
same $a$-type weight to any tiling. Without loss of generality, we can
therefore henceforth set $a_{o(t)}^{w(t)} = 1$ for any $t$ such that
$w(t) \in {\cal B}$.

It is convenient to illustrate the distinction of Fig.~\ref{fig3} in
another way, by decorating the interior of each tile by an oriented
half-arch (see Fig.~\ref{fig4}). We then have the following \\

\begin{lemma}
In any valid tiling, the orientation of the half-arches is conserved across the junctions of the tiles.
\end{lemma}

\begin{proof}
  By rotational and reflectional symmetry, it suffices to consider the
  junction between a tile in the $a_1$-type position
  (cf.~Fig.~\ref{fig1}) and the neighbouring tile immediately
  North-East of it.

\begin{figure}
  \centering
  \includegraphics[width=200pt]{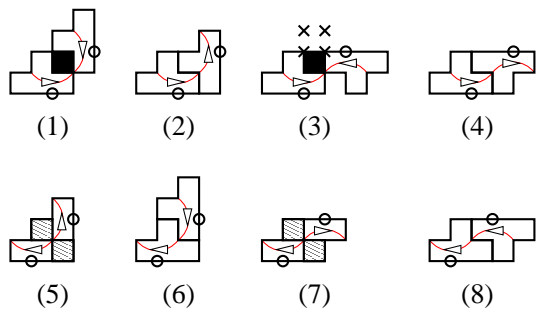}
  \caption{Arrow conservation at tile junctions.}
  \label{fig5}
\end{figure}

Consider then the situation where the first tile is of the $b_1$ type
(cf.~Fig.~\ref{fig4}).  The position of the neighbouring tile is
constrained by Theorem 1 to have the correct location of its white
vertex. When the neighbouring tile has its long side vertical, there
are just two possibilities, shown in Fig.~\ref{fig5}, panels (1) and
(2). The possibility (1) is disallowed since it leaves an untiled hole
(shown in black), whereas (2) is allowed.  When the neighbour tile's
long side is horizontal, we have the possibilities (3) and (4). But
(3) is actually disallowed, since, when adding a further tile to cover
the lattice face shown in black, a white point would be placed in one
of the four vertices marked by a cross, and neither of these vertices
is a valid position for a white point according to Theorem 1. In
summary, only the possibilities (2) and (4) are allowed, and both of
those conserve the arrow orientation across the tile junction.

Consider next the situation where the first tile is of the $b_2$ type.
A neighbour tile with a vertical long side leads to the possibilities
(5) and (6), but (5) is disallowed since it leads to two overlaps
(shown hatched). Similarly, when the neighbour's long side is horizontal,
one can rule out (7) due to two overlaps, whereas (8) is allowed. In
summary, only (6) and (8) are allowed and compatible with arrow
conservation.
\end{proof}

\begin{defn}
Let $G=(V,E)$ be a finite undirected graph with vertex set $V$ and edge
set $E$. The multivariate Tutte polynomial of $G$ with parameter
$Q$ and edge weights ${\bf v}=\{v_e\}_{e\in E}$ is
\begin{equation}
 Z_G(Q,{\bf v}) = \sum_{A \subseteq E} Q^{k(A)} \prod_{e \in A} v_e \,,
 \label{Tutte}
\end{equation}
where $k(A)$ denotes the number of connected components in the subgraph
$(V,A)$.
\end{defn}

For the relation of the multivariate Tutte polynomial to the more standard
two-parameter Tutte polynomial $T_G(x,y)$ we refer to Section 2.5 of
\cite{Sokal}. In the physics literature, $Z_G(Q,{\bf v})$ is better
known as the partition function of the $Q$-state Potts model in the
Fortuin-Kasteleyn formulation \cite{FK}.

In the remainder of the paper we let $V = B_{\rm even}$ and $E$ the associated
natural set of horizontal and vertical edges of length $4$. In other words, $G$
is an $(m-1) \times (n-1)$ rectangle $G$, as referred to in the abstract.
Furthermore, let $w(e)$ be the natural bijection between edges $e \in E$
and white vertices $w \in W \setminus {\cal B}$ that are {\em not} on the
boundary ${\cal B}$ of ${\cal D}$.

We can now state our main result which is

\begin{theor}
The generating function $F_{\cal D}(\{a\},\{b\})$ of T-tetromino tilings of
the $4m \times 4n$ rectangle ${\cal D}$ depends only on the parameters
$\{a\}$ through the combinations $a_1 a_3$ and $a_2 a_4$.
For the specific choice of tile weights
\begin{equation}
 a_1 a_3 = \left \lbrace \begin{array}{ll}
 x_e & \mbox{{\rm for} $w(e) \in W_{\rm even} \setminus {\cal B}$} \\
 1   & \mbox{{\rm otherwise}}
 \end{array} \right.
 \qquad
 a_2 a_4 = \left \lbrace \begin{array}{ll}
 x_e & \mbox{{\rm for} $w(e) \in W_{\rm odd} \setminus {\cal B}$} \\
 1   & \mbox{{\rm otherwise}}
 \end{array} \right.
 \label{aweights}
\end{equation}
and
\begin{equation}
 b_1 = (b_2)^{-1} = q^{1/4}
 \label{bweights}
\end{equation}
we have the identity
\begin{equation}
 Q^{m n/2} F_{\cal D}(\{a\},\{b\}) = Z_G(Q;{\bf v})
 \label{main_result}
\end{equation}
with the correspondence between parameters
\begin{equation}
 Q = (q+q^{-1})^2, \qquad
 v_e = (q+q^{-1}) x_e \mbox{ {\rm for} $e \in E$}
\end{equation}
\end{theor}

\begin{proof}
By Theorem 1, any white vertex $w \in W \setminus {\cal B}$ is at the junction
between the long sides of exactly two tiles. If the long side is horizontal,
the contribution of the $a$-type weight is $a_1^w a_3^w$, and if it is
vertical the contribution is $a_2^w a_4^w$. This proves the first part of
the theorem. (Above we have already explained that one may set all
$a_i^w = 1$ when $w \in W \cap {\cal B}$.)

\begin{figure}
  \centering
  \includegraphics[width=120pt]{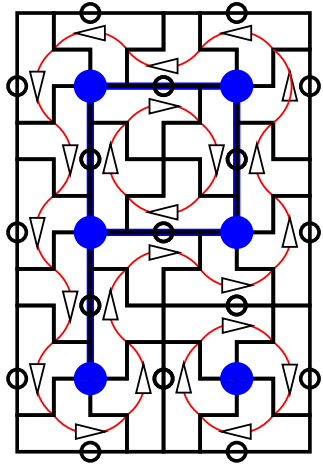}
  \caption{A valid tetromino tiling of the domain ${\cal D}$ shown in
Fig.~\ref{fig1}. The vertices $V$ entering the definition of the Tutte
polynomial (\ref{Tutte}) are shown as fat solid circles, and the edges
of the subset $A \subseteq E$ are shown as fat lines.}
  \label{fig6}
\end{figure}

By Lemma 1 the oriented half-arches form continuous curves with a
consistent orientation (see Fig.~\ref{fig6} for an example). Moreover,
by the consistency of the position of the white vertices
(cf.~Figs.~\ref{fig2}--\ref{fig3}), these curves cannot end on ${\cal
B}$. Therefore, the half-arches form a set of oriented cycles.

Let $C$ be the set of all possible configurations of cycles arising
from valid tilings of ${\cal D}$, but disregarding the orientations of
the cycles. The summation in the generating function (\ref{genF}) can
then be written $\sum_{{\cal T}({\cal D})} = \sum_{C} \sum_{{\cal
    T}({\cal D})|C}$, where the last sum is over cycle orientations only. A
close inspection of panels (2), (4), (6) and (8) in Fig.~\ref{fig5}
reveals that the orientation of a given oriented cycle $c_0$ can be
changed independently of all other oriented cycles in the tiling $T
\in {\cal T}({\cal D})$ by shifting the tiles traversed by $c_0$ by
one lattice unit to the left, right, up or down, but without moving
any tiles not traversed by $c_0$. In particular, this transformation
does not change the orientation of any tile (cf.~Fig.~\ref{fig1}), and
so does not change the $a$-type weights. We have therefore
\begin{equation}
 F_{\cal D}(\{a\},\{b\}) = \sum_{C} \prod a_{o(t)}^{w(t)}
 \sum_{{\cal T}({\cal D})|C} \prod_{j=1}^2 b_j^{B_j} \,.
 \label{Fgen1}
\end{equation}

[The remainder of the proof parallels a construction used in \cite{BKW}
to show the equivalence between the Potts and six-vertex model partition
functions.]

The partial summation $\sum_{{\cal T}({\cal D})|C} \prod_{j=1}^2
b_j^{B_j}$ appearing in (\ref{Fgen1}) can be carried out by noting
that the $b$-type weights (\ref{bweights}) simply amount to weighting
each turn of $c_0$ through an angle $\pm \pi/2$ by the factor $q^{\pm
1/4}$. Since the complete cycle turns a total angle of $\pm 2\pi$
depending on its orientation, the
product of the $b$-type weights associated with the cycle $c_0$,
summed over its two possible orientations, yields the weight $q+q^{-1}
= Q^{1/2}$.  Denoting by $\ell(C)$ the number of cycles in $C$, we
have therefore found that
\begin{equation}
 \sum_{{\cal T}({\cal D})|C} \prod_{j=1}^2 b_j^{B_j} = Q^{\ell/2} \,.
\end{equation}

It now remains to sum over the un-oriented cycles $C$. There is an
easy bijection between $C$ and the edge subsets $A \subseteq E$
appearing in (\ref{Tutte}). Namely, an edge $e\in E$ belongs to $A$
iff it is not cut by any cycle in $C$ (see
Fig.~\ref{fig6}). Conversely, given $A$, the cycles are constructed so
that they cut no edge in $A$, cut all edges in $E \setminus A$, and
are reflected off the boundary ${\cal B}$ in the vertices $W \cap {\cal
B}$.  Under this construction, the $a$-type weights (\ref{aweights})
simply give $\prod_{e \in A} x_e$. Collecting the pieces this amounts
to
\begin{equation}
 F_{\cal D}(\{a\},\{b\}) = \sum_{A \subseteq E}
 Q^{\ell(A)/2} \prod_{e \in A} \frac{v_e}{Q^{1/2}} \,.
\end{equation}

Using finally the Euler relation, and remarking that $|V|= m n$,
one obtains (\ref{main_result}).
\end{proof}

Let us discuss a couple of special cases of Theorem 2. First, when
there are no $a$-type weights (i.e., setting $a_{o(t)}^{w(t)} = 1$
for any $o(t)$ and $w(t)$), the Tutte polynomial has $v_e = q+q^{-1}$ for
any $e \in E$. The tiling entropy defined by
\begin{equation}
 S_G(q) \equiv \lim_{m,n \to \infty} \big( F_{\cal D} \big)^{\frac{1}{m n}}
\end{equation}
is then related to the one found for the Tutte polynomial ({\em alias}
Potts model) by Baxter \cite{Baxter73} using the Bethe Ansatz
method. For $0 < Q < 4$ one sets $q \equiv {\rm e}^{i \mu}$ with $\mu
\in (0,\pi/2)$ and finds
\begin{equation}
 \log S_G(q) = \int_{-\infty}^{\infty} {\rm d}t \,
 \frac{\sinh\left[(\pi-\mu)t \right] \, \tanh(\mu t)}
      {t \, \sinh(\pi t)} \,,
\end{equation}
whereas for $Q > 4$ one sets $q \equiv {\rm e}^\lambda$ and finds
\begin{equation}
 \log S_G(q) = \lambda +
 2 \sum_{n=1}^{\infty} \frac{{\rm e}^{-n \lambda} \, \tanh(n \lambda)}{n} \,.
\end{equation}
Finally, for $Q=4$ one has
\begin{equation}
 S_G(1) = \left( \frac{\Gamma(5/4)}{\Gamma(3/4)} \right)^4 \,.
\end{equation}

An even more special case of the identity (\ref{main_result}) arises when
$b_1=b_2=1$ as well. $F_{\cal D}$ is then the unweighted sum over
all T-tetromino tilings of ${\cal D}$. In terms of the standard
Tutte polynomial $T_G(x,y)$ one then has $F_{\cal D} = 2 T_G(3,3)$,
as first proved by Korn and Pak \cite{KornPak}.


\begin{thebibliography}{99}

\bibitem{Walkup} D.W. Walkup,
 Amer.~Math.~Monthly {\bf 72}, 986--988 (1965).

\bibitem{Sokal} A.D. Sokal, in {\em Surveys in Combinatorics, 2005},
 edited by Bridget S. Webb (Cambridge University Press, 2005), pp.~173--226;
 math.CO/0503607.

\bibitem{FK} C.M. Fortuin and P.W. Kasteleyn,
 Physica {\bf 57}, 536--564 (1972).

\bibitem{BKW} R.J. Baxter, S.B. Kelland and F.Y. Wu,
 J. Phys. A {\bf 9}, 397--406 (1976).

\bibitem{Baxter73} R.J. Baxter,
 J. Phys. C: Solid State Phys. {\bf 6} L445--L448 (1973).

\bibitem{KornPak} M. Korn and I. Pak,
 Theor.~Comp.~Science {\bf 319}, 3--27 (2004).

\end{thebibliography}
\end{document}